\documentclass{amsart}
\usepackage{amsmath,amssymb, amsbsy}
\usepackage[dvips]{graphicx}
\usepackage{color}
\usepackage[latin1]{inputenc}
\usepackage[active]{srcltx}
\usepackage{wrapfig}
\usepackage{graphicx}
\usepackage[usenames,dvipsnames]{pstricks}
\usepackage{epsfig}
 \usepackage{pst-grad} 
 \usepackage{pst-plot} 
\usepackage{epic}

\theoremstyle{plain}

\newtheorem{remark}{Remark}

\numberwithin{equation}{section}
\long\def\salta#1{\relax}
\newcommand{\R}{{I\!\!R}}

\newcommand{\p}{\partial}

\newcommand{\re}{{I\!\!R}}
\newcommand{\ren}{\re^N}

\newcommand{\dyle}{\displaystyle}

\newcommand{\io}{\int\limits_\O}

\newcommand{\sob}{W_0^{1,2}}

\renewcommand{\a }{\alpha }
\renewcommand{\b }{\beta }

\newcommand{\D }{\Delta }
\newcommand{\e }{\varepsilon }

\newcommand{\g }{\gamma}

\renewcommand{\l }{\lambda }
\renewcommand{\L }{\Lambda }

\newcommand{\n }{\nabla }

\newcommand{\s }{\sigma }

\renewcommand{\O }{\Omega }

\newtheorem{Theorem}{Theorem}[section]

\newtheorem{Definition}[Theorem]{Definition}

\newcommand{\cqd}{{\unskip\nobreak\hfil\penalty50
        \hskip2em\hbox{}\nobreak\hfil\mbox{\rule{1ex}{1ex} \qquad}
        \parfillskip=0pt \finalhyphendemerits=0\par\medskip}}

\begin{document}


\title[Existence of positive solutions to a nonlinear elliptic system]{Existence of positive solutions to a nonlinear\\ elliptic system with nonlinearity involving gradient term }

\author[B. Abdellaoui, A. Attar, EL-HAJ LAAMRI ]{Boumediene Abdellaoui, Ahmed Attar, El-Haj Laamri}
\thanks{ The first author is partially supported by project MTM2016-80474-P, MINECO, Spain. } 

\address{\hbox{\parbox{5.7in}{\medskip\noindent {B. Abdellaoui, A. Attar, Laboratoire d'Analyse Nonlin\'eaire et Math\'ematiques
Appliqu\'ees. \hfill \break\indent D\'epartement de
Math\'ematiques, Universit\'e Abou Bakr Belka\"{\i}d, Tlemcen,
\hfill\break\indent Tlemcen 13000, Algeria.}}}}
\address{\hbox{\parbox{5.7in}{\medskip\noindent{EL-HAJ LAAMRI, Institut Elie Cartan,\\ Universit\'{e} Lorraine,\\ B. P. 239, 54506 Vand{\oe}uvre l\'es
 Nancy, France. \\[3pt]
        \em{E-mail addresses: }\\{\tt boumediene.abdellaoui@inv.uam.es, \tt ahm.attar@yahoo.fr, \tt el-haj.laamri@univ-lorraine.fr}.}}}}

\date{\today}

\begin{abstract}
In this work we analyze the existence of solutions to the nonlinear elliptic system:
\begin{equation*}
\left\{
\begin{array}{rcll}
-\Delta u & = & v^q+\a g & \text{in }\Omega , \\
-\Delta v& = &|\nabla u|^{p}+\l f &\text{in }\Omega , \\
u=v&=& 0 & \text{on }\partial \Omega ,\\
u,v& \geq & 0 & \text{in }\Omega,
\end{array}%
\right.
\end{equation*}
where $\Omega$  is a bounded domain of $\ren$ and $p\ge 1$, $q>0$ with $pq>1$. $f,g$ are nonnegative measurable functions with additional hypotheses and $\a, \l\ge 0$.

As a consequence we show that the fourth order problem
\begin{equation*}
\left\{
\begin{array}{rcll}
\Delta^2 u & = &|\nabla u|^{p}+\tilde{\l} \tilde{f} &\text{in }\Omega , \\
u=\D u&=& 0 & \text{on }\partial \Omega ,\\
\end{array}%
\right.
\end{equation*}
has a solution for all $p>1$, under suitable conditions on $\tilde{f}$ and $\tilde{\l}$.

\end{abstract}

\maketitle

\section{Introduction}\label{section00}

The aim of this paper is to discuss the existence of solutions to the following elliptic system
\begin{equation}\label{S}
\left\{
\begin{array}{rcll}
-\Delta u & = & v^q+\a g & \text{in }\Omega , \\
-\Delta v& = &|\nabla u|^{p}+\l f &\text{in }\Omega , \\
u=v&=& 0 & \text{on }\partial \Omega ,\\
u,v & \geq& 0 & \text{in }\Omega,
\end{array}%
\right.
\end{equation}
where $\Omega\subset \ren$ is a bounded domain. We will consider the case $p\ge 1, q>0$ with $pq>1$,
$f,g$ are nonnegative measurable functions and $\a, \l$ are nonnegative real constants.

Our goal is to get \textit{natural} conditions on the parameter $\a, \l$ and the data $f, g$, in  order to prove the existence of positive solutions to the problem \eqref{S} under the condition $pq>1$. By solution, we mean solution  in the sense of distributions (see definition \ref{def-faible}).
\vskip3mm
The class of elliptic systems with gradient term appears when considering electrochemical models in engineering and some other model in fluid dynamics. We refer to \cite{C}, and \cite{Diez} for more details and more applications of this class of systems.
\vskip3mm
Existence results for nonlinear elliptic systems with gradient term are well known in some particular cases. For example, in the case where  $\alpha=\lambda=1$, the authors in \cite{AtB} established  that system (\ref{S}) has a solution for all $(f,g)\in L^2(\Omega)\times L^2(\Omega)$
under the condition that $pq<1$ and $0<p<2$.
In \cite{BOPorr}, Boccardo-Orsina-Porretta  investigated  the system
\begin{equation*}
\left\{
\begin{array}{rcll}
- \text{div}(b(x,z)\n u)&=& f(x) &\text{ in }\Omega , \\
- \text{div}(a(x,z)\n z)&=&b(x,z)|\nabla u|^{2} &\text{ in }\Omega , \\
u=z&=&0  &\text{ on }\partial \Omega,
\end{array}%
\right.
\end{equation*}
where the  functions $(x,s)\mapsto a(x,s), b(x,s)$ are positive and coercive Carath\'eodory functions. Under the hypothesis that $f\in L^m(\O)$, with $m\ge \frac{2N}{N
+2}$, they proved the existence and the regularity of a positive solution. In  \cite{BOPu}, Boccardo-Orsina-Puel studied  the  system
\begin{equation*}
\left\{
\begin{array}{rcll}
-\text{div}(a(x,z)\n u)&=& f(x) & \text{ in }\Omega , \\
-\text{div}(b(x)\n z)+K(x,z)|\nabla u|^{2}&=&g(x) & \text{ in }\Omega , \\
z=u&=&0 & \text{  on }\partial \Omega,
\end{array}%
\right.
\end{equation*}
 where the gradient appears as an absorption term.
\vskip3mm
 \noindent It is clear that for $g=0$, setting $v=|-\D u|^{\frac 1q-1}(-\D u)$, the system \eqref{S} is reduced to the following fourth order problem
 \begin{equation}\label{PBpp}
\left\{
\begin{array}{rcll}
-\Delta(|-\D u|^{\frac 1q-1}(-\D u))&= & |\nabla u|^{p}+\l f &
\text{ in }\Omega , \\ u=\D u &=& 0 &\hbox{  on } \p\Omega,\\ u&>&0 &\hbox{
	in }\Omega.
\end{array}%
\right.
\end{equation}
If $q=1$, the above problem is reduced to the following one
\begin{equation}\label{PBi}
\left\{
\begin{array}{rcll}
\Delta^2 u &= & |\nabla u|^{p}+\l f &
\text{ in }\Omega , \\ u=\D u &=& 0 &\hbox{  on } \p\Omega,\\ u&\geq &0 &\hbox{
	in }\Omega.
\end{array}%
\right.
\end{equation}
Hence, to get a positive solution to problem \eqref{PBi}, we  just have to show that the system (\ref{S}) has a solution with $q=1$ and $\a=0$.\\
Problems related to the Bi-Laplacian operator are widely studied in the literature, we refer to \cite{GG} and the references therein. We refer also to the paper \cite{EP} where an extension of the Kardar-Parisi-Zhang equation for the Bi-Laplacian operator is also studied.
\vskip2mm
Notice that if we substitute the fourth order operator $\D^2 u$ by the classical Laplacian, the previous problem takes the form
\begin{equation}\label{lapl}
\left\{
\begin{array}{rcll}
-\Delta u&=&|\nabla u|^{p}+\l f &\hbox{ in }\Omega , \\ u&=&0 &\hbox{  on } \p\Omega.
\end{array}%
\right.
\end{equation}
This problem has been widely studied in the literature, we refer to \cite{LL, BMPu, BGM,ChoCho, FM1, BarlPorre, GreTrom, PLL} and references therein.
%
If $p\le 2$, using suitable comparison principle and under suitable hypothesis on $f$, Alaa-Pierre  proved  in \cite{Alp} the existence of positive solutions to \eqref{lapl} that is in a suitable Sobolev space. Some extension were proved in \cite{GMP} using truncations arguments. The existence result in \cite{GMP} holds for a more general class of elliptic operators. However, in the case of the Laplacian operator, the condition $p\le 2$ seems to be optimal and cannot be improved by using the techniques of \cite{GMP}. 
 On the other hand, if $u\in \sob(\O)$, then in general, $\D(T_k(u))\notin L^1(\O)$ and it seems to be more complicated to handle  problem \eqref{PBi} using truncation arguments.
The existence result obtained in \cite{HMV} is more general because it holds for all $p>1$. However, it seems that the techniques of the proofs are strongly related to a fine estimate on the Green function and cannot be extended to more general operators.\\
In the case where $\O=B_R(0)$, then for $p=1, q=2$ and $g=f=0$, based on ODEs technics, the authors in \cite{Diez} were able to show that the corresponding system has a unique radial large solution.\\

Our paper mainly complements the investigations of \cite{AtB}. However as it was noticed in \cite{AtB}, the arguments used to deal with the case $pq<1$ can not be adapted to the case $pq>1$. Some new arguments are needed to deal with the new situation. Hence, we will use in a convenient way, the Schauder fixed point theorem, following some ideas from Phuc in \cite{Phuc}. Moreover, in order to show that (\ref{S}) has a solution, it is necessary  to have regularity assumptions on the data $f,g$ and smallness conditions on $\lambda$ and $\alpha$. Indeed, such conditions are necessary even in the case of a single  equation. To help understand the situation, let us consider the following semilinear equation :
\begin{equation*}
(P_r)\left\{
\begin{array}{rcll}
-\Delta w & = & w^r+ h & \text{in }\Omega , \\
w&=& 0 & \text{on }\partial \Omega ,\\
w& \geq & 0 & \text{in }\Omega,
\end{array}
\right.
\end{equation*}
where $r\in (1,+\infty)$ and $h\geq 0$. As it is proved by Baras-Pierre in \cite{BPihp}, two conditions on $h$ are necessary
for the existence of solutions to $(P_r)$ :\\
(i) A regularity condition : $h$ should  be ``regular" enough. \\
(ii) A size condition : even if $h\in \mathcal{C}_0^\infty (\Omega)$, it should be small enough. For instance, if $h=\gamma\psi$ where $\psi\in \mathcal{C}^\infty_0(\O)$, $\psi\gneqq 0$ and $\gamma >0$,  then there exists $\gamma^*$ so that $(P_r)$ does not have any solution for $\gamma>\gamma^*$.\vskip2mm
\noindent In other words, a necessary and sufficient condition on $h$ for the existence of solutions to $(P_r)$ can be formulate in a simple way, saying that 
a certain  ``norm" of $h$  should be exactly less than or equal to
 $k(r)=\displaystyle\frac{r-1}{r^{r'}}$ where $ r'=\displaystyle\frac{r}{r-1}$. This quantity is defined by duality through a functional which is naturally associated with the above problem.
  Mathematically speaking, one has the following result :
\begin{Theorem}[\textbf{Baras-Pierre\cite{BPihp}}]
 Let $h$ be a nonnegative measurable function on $\Omega$. Then
$(P_r)$ \textit{ has a weak solution if and only if }
 \begin{equation*}
(H)\left\{
\begin{array}{lll}
\forall \varphi \in \mathcal{C}^{\infty}(\Omega) \text{ with } \varphi\gneqq 0 \text{ in } \Omega \text{ and }  \varphi = 0 \text{ on } \partial\Omega ,\\
\displaystyle\int_\Omega h\varphi\leq k(r) \int_\Omega |\Delta\varphi|^{r'} |\varphi|^{1-r'}.
\end{array}
\right.
\end{equation*}
\end{Theorem}

Thus, we will deal with this type of functionals in order  to describe exactly the optimal size of the data $f$, $g$ and to determine the largest set of $(\lambda,\alpha)$ such that the system $(\ref{S}) $ has a solution.
\vskip4mm
The paper is divided into four sections. We  state our main results in section \ref{sec-2}. Then we introduce some  useful tools in section \ref{sec-3}. Finally  the last section  is devoted to prove our theorems and we also give  some nonexistence results that show, in some sense, the optimality of the hypothesis imposed in Theorem \ref{bih}.
\newpage

\section{Main results}\label{sec-2}
Throughout this paper, $\Omega$ is a bounded open set of $\R^N$ and $f,g :\Omega \rightarrow \R$ are nonnegative measurable functions.
\vskip2mm
Let us precise what we mean by a solution.
\begin{Definition}\label{def-faible}
Assume that $p\ge 1$ and $q>0$. Let $f, g$ be nonnegative functions such that $f,g\in L^1(\O)$. We say that $(u,v)\in L^1(\O)\times L^1(\O)$, with $u,v\ge 0$, is a weak solution to system \eqref{S} if $(u,v)\in W^{1,p}_0(\Omega)\times W^{1,1}_0(\O), v^q \in L^1(\O)$ and for all
$\varphi, \psi\in \mathcal{C}^\infty_0(\O)$, we have
\begin{equation}\label{eq1f}
\dyle\io \nabla \varphi \nabla u
=\dyle\io v^q\varphi+\a\dyle\io g \varphi,\mbox{  and   }\dyle\io \nabla\psi\nabla v=\dyle\io| \nabla u| ^p\psi+\l\dyle\io f\psi.
\end{equation}
\end{Definition}
%

We are now ready to state our main existence results.

\begin{Theorem}\label{syst1}
Assume that $p\ge 1$, $q>0$ with $pq>1$. Let $(m,\s)\in (1,+\infty)^2$. Suppose that $(f,g) \in L^{m}(\O)\times L^{\s}(\Omega)$ where  $(m,\s)$ satisfies one of the following  conditions
\begin{equation}\label{condi0}
\left\{
\begin{array}{rcll}
m,\s &\in & (1, N),\\
pm & < & \dfrac{\s N}{N-\s}=\s^*,\\
\dfrac{q\s}{N+q\s} &<& \dfrac{m}{N-m},
\end{array}%
\right.
\end{equation}
or
\begin{equation}\label{condi01}
m\ge N\mbox{  and  }\s>\dfrac{pmN}{N+pm},
\end{equation}
or
\begin{equation}\label{condi03}
\s \ge N\mbox{  and  }m>\dfrac{q\s N}{N+2q\s}.
\end{equation}
Then there exists $\Lambda^*>0$ such that for all $(\lambda, \alpha)\in \Pi$ where
\begin{equation}\label{pi1}
\Pi : =\{(\lambda, \alpha) \in [0,+\infty)\times [0,+\infty) \;|\; \lambda\|f\|_m+\alpha^p \|g\|_\s \leq \L^*\},
\end{equation}
 the system (\ref{S}) has a nonnegative solution $(u,v)$. Moreover $(u,v)\in W^{1,\theta}_0(\O)\times W^{1, r}_0(\O)$ for all $\theta<\dfrac{\s N}{(N-\s)_+}$ and $r<\dfrac{mN}{(N-m)_+}$.
\end{Theorem}
\vskip2mm
\begin{remark} To give some light on the hypothesis (\ref{condi0}), (\ref{condi01}) and (\ref{condi03}), let us explicit the size conditions on $(p,q)$ for a given $(m,\s)$.\\
If  $m=\s=2$ and  $N\geq 3$,  then  condition (\ref{condi0}) is satisfied for $p<\dfrac{N}{N-2}$ and  $q<\dfrac{N}{(N-4)_+}$ ;\\
If  $m=\s=N/2$ and $N\geq 3$,  then  condition (\ref{condi01}) is satisfied for $p<2$ and  for all $q$.\\
\noindent If  $m=\s\geq N$, then the  conditions (\ref{condi01}) and (\ref{condi03}) are satisfied for all $(p,q)\in [1,+\infty)^2$.\vskip1mm
\noindent If $\s=N$ and $m=N/2$, then the  conditions (\ref{condi01}) and (\ref{condi03}) are satisfied for all $(p,q)\in [1,+\infty)^2$.\vskip1mm
\noindent If $m=N$ and $\s=\dfrac{3N}{4}$, then the  condition (\ref{condi01}) is satisfied for all $p<3$ and all $q$.
\end{remark}


Notice that the set $\Pi$ defined in \eqref{pi1} is a bounded set of $\re^2_+$. The next nonexistence result explain clearly that  a smallness condition on  $(\lambda, \alpha)$ is necessary for the existence of solutions to (\ref{S}), at least for $p>1$ and $q\ge 1$.

\begin{Theorem}\label{smallness_Lambda-Alpha} 
Suppose that $p>1$ and $q\ge 1$. Let $f,g$ be  nonnegative measurable functions such that $(f,g)\neq (0,0)$  and $(\lambda, \alpha)\in (0,+\infty)^2$.
Assume that system $(\ref{S})$ has a nonnegative solution $(u,v)$, then there exists $(\lambda^*, \alpha^*)\in (0,+\infty)^2$ such that $\lambda \leq \lambda^*$ and $\alpha \leq \alpha^*$.
\end{Theorem}
%
%
In the particular case where $q=1$ and $\a\equiv 0$, and as a direct application of the Theorem  \ref{syst1}, we obtain the following existence result
for the Bi-Laplacian problem with gradient term.

\begin{Theorem}\label{bih}
Let $p>1$. Suppose $0\lneqq f\in L^m(\O)$ where  $m >\max\{1,\frac{N}{3p'}\}$. Then there exists $\l^*>0$ such that if $\l<\l^*$ the problem 
\begin{equation}\label{PBih}
\left\{
\begin{array}{rcll}
\Delta^2 u &= & |\nabla u|^{p}+\l f &
\text{ in }\Omega , \\ u=\D u &=& 0 &\hbox{  on } \p\Omega,\\ u&\geq &0 &\hbox{
	in }\Omega.
\end{array}%
\right.
\end{equation}
has a solution $u$ such that $u\in W^{4,m}(\O)$ if $m<\frac{N}{2}$ and $u\in \mathcal{C}^{s}(\O)$ with $s<4-\frac{N}{m}$ if $m>\frac N2$.
\end{Theorem}
Finally we give  a nonexistence result that, in some sense, justifies
the regularity conditions imposed on $f$ to get the existence of nonnegative solution to problem \eqref{PBih}.

\begin{Theorem}\label{Non100}
Assume that  $1\leq m<\max\{1,\frac{N}{3p'}\})$ where $p>1$, then there exists $f\in L^m(\O)$ with $f\gneq 0$ such that problem \eqref{PBih} has non positive solution for all $\lambda>0$.
\end{Theorem}

%
%
%
%
\section{ Useful results}\label{sec-3}
%
For the convenience of the reader and for the sake of completeness, we recall in this short section some classical  results we will use in our proofs.
\vskip2mm
In order to prove the main existence result of this paper we will use useful theorems: Schauder fixed point Theorem and Vitali's Theorem.

\begin{Theorem}[\textbf{Schauder fixed point theorem}] \label{pointfixe}
Assume that $E$ is a closed convex set of a Banach space $X$. Let $T$ be a continuous and compact mapping from $E$ into itself. Then $T$ has a fixed point in $E$.
\end{Theorem}
%

\begin{Theorem}{\bf (Vitali)}\label{vitali} 
Let $(E,\mu)$ be a measured space such that $\mu(E)<+\infty$, let  $1\leq p<+\infty$  and let $\{f_n\}_n\subset L^p(E)$  such that $f_n\rightarrow f\;\; a.e$. If $\{f_n^p\}_n$ is uniformly integrable over $E$, then $f\in L^p(E)$ and $f_n\rightarrow f$ in $L^p(E)$.
\end{Theorem}

Systematically , we will use the following regularity result proved in \cite[Appendix]{BPaif}.
\begin{Theorem}\label{th0} 
Assume that $h\in L^1(\O)$, then the problem
\begin{equation}
\left\{
\begin{array}{rcll}
-\Delta z &=& h & \mbox{ in }\Omega, \\ z &=& 0 & \mbox{  on } \p\Omega
\end{array}%
\right.  \label{PQ}
\end{equation}
has a unique weak  solution $z\in W^{1,s}_0(\Omega)$ for all $s\in [1, \frac{N}{N-1})$.\\ Moreover, for $s\in [1, \frac{N}{N-1})$ fixed, there exists a positive constant $C=C(\Omega, N,s)$ such that
\begin{equation}\label{dd22}
||\n z||_{L^s(\Omega)}\le C||h||_{L^1(\O)}. 
\end{equation}
and the operator $\Phi : h\mapsto z$ is compact from $L^1(\O)$ into $W^{1,s}_0(\Omega)$.
\end{Theorem}

Finally, let us recall the following classical regularity result that we  will use in several proofs below.

\begin{Theorem}\label{key-2} 
Let $h\in L^m(\O)$ with $m> 1$. Then the problem
\begin{equation}\label{gener1}
\left\{
\begin{array}{rcll}
-\Delta z &= & h &
\mbox{ in }\Omega , \\ z &=& 0 &\hbox{  on } \p\Omega.
\end{array}%
\right.
\end{equation}
has a unique weak solution $z$. Moreover there exists a positive constant $C =C(\Omega, N,m)$ independent of $h$ such that
\begin{enumerate}
\item If $1<m<N$, then
\begin{equation}\label{dd11}
||\n z||_{L^{m^*}(\Omega)}\le C||h||_{L^m(\O)}\mbox{  where  }m^*=\frac{m N}{N-m}.
\end{equation}

\item If $m=N$,  $|\n z|\in L^s(\Omega)$ for all $s\in [1,+\infty)$.
\item If $m>N$,  $z\in \mathcal{C}^{1, \gamma}(\Omega)$ for some $\gamma\in (0,1)$.
\end{enumerate}
\end{Theorem}
For a proof see for instance \cite{Bre} or \cite{BocCroce}.

\begin{remark}
As a consequence of the above theorems, we can prove that for all $h\in L^1(\O)$, there exists a unique solution $w$ of the problem
\begin{equation}\label{ggent}
\left\{
\begin{array}{rcll}
\Delta^2 w &= & h &
\text{ in }\Omega , \\ w=\D w &=& 0 &\hbox{  on } \p\Omega,
\end{array}%
\right.
\end{equation}
with $w\in W^{1,\theta}_0(\O)$ for all $\theta<\frac{N}{N-3}$ and
$$
||w||_{W^{1,\theta}_0(\O)}\le C(\O, \theta)||h||_{L^1(\O)}.
$$
\end{remark}

\section{Proofs of the main results. }\label{sec-4}

\subsection{Proof of theorem \ref {syst1}}

We give the proof in the case where $(m, \s)$ satisfying \eqref{condi0} i.e.
$$m,\s \in (1, N),\;
pm <\dfrac{\s N}{N-\s}=\s^*,\;
\dfrac{q\s}{N+q\s} <\dfrac{m}{N-m}.$$

 The other cases follow by using the same arguments.

For $s\ge 0$, we define the function
$$
\Upsilon(s)=s^{\frac{1}{pq}} -\tilde{C}s,
$$
where $\tilde{C}$ is a universal positive constant that depends only on datum  and that  will specify later.\\
Using the  fact that $pq>1$, then there exists $s_0>0$ such that $\Upsilon(s_0) =0$, $\Upsilon(s) >0, \; \forall s\in (0,s_0)$, $\Upsilon(s) <0,\; \forall s\in (s_0, +\infty)$  and we get the existence of a positive constants   $\ell>0$ and $\L^*>0$  such that
$$
\max_{s\ge 0}\Upsilon(s)=\Upsilon(\ell)=\L^*.
$$
Thus
$$
\ell^{\frac{1}{pq}}=\tilde{C}(\ell+\frac{\L^*}{\tilde{C}}).
$$
Fix $\ell>0$ as above and define the set
$$
\Pi\equiv \bigg\{(\l, \a)\in [0,+\infty)\times [0,+\infty) \;;\; \l||f||_{L^{m}(\O)}+\a^p||g||_{L^{\s}(\O)}^p\le \frac{\L^*}{\tilde{C}}\bigg\}.
$$
\vskip2mm
\noindent It is clear that $\Pi$ is non empty, bounded  and for all $(\l,\a)\in \Pi$, we have
\begin{equation}\label{relation-entre-l-f-g-lambda, alpha}
\tilde{C}\bigg(\ell+\l||f||_{L^{m}(\O)}+\a^p||g||_{L^{\s}(\O)}^p\bigg)\le \ell^{\frac{1}{pq}}.
\end{equation}
\vskip2mm
\noindent
In what follows we fix $(\l,\a)\in \Pi$.
Since $\frac{q\s N}{N+q\s}<\frac{mN}{N-m}$, there exists $r>1$ such that
\begin{equation}\label{mad}\frac{q\s N}{N+q\s}<r<\frac{mN}{N-m}.
\end{equation}
Moreover, if $r<N$, then
\begin{equation}\label{TTT}
\s q<\frac{rN}{N-r}=r^*.
\end{equation}
and if $r\ge N$, then \eqref{TTT} holds trivially with $r^*=\infty$.\vskip3mm
\noindent Now let us fix $\ell$ and $r$ as above, we define the set
\begin{equation}\label{sett}
E=\{w\in W^{1,1}_0(\O)\; ;\; w\in W^{1,r}_0(\O)\mbox{  and  }||\n w||_{L^{r}(\Omega)}\le \ell^{\frac{1}{pq}}\}.
\end{equation}
One can easily verify that $E$ is a closed convex subset of $W^{1,1}_0(\O)$. 
Let us consider the operator
$$T:  E \longrightarrow W^{1,1}_0(\O)\  \  $$
$$w \longmapsto T(w)=v$$
where $v$ is the unique solution to problem
\begin{equation}\label{P100}
\left\{
\begin{array}{rcll}
-\Delta v & = & |\nabla u|^{p}+\l f
& \text{ in }\Omega,\\
v &=& 0  & \text{ on }\partial \Omega,
\end{array}
\right.
\end{equation}
with $u$ being the unique solution to the  problem
\begin{equation}\label{P110}
\left\{
\begin{array}{rcll}
-\Delta u &= & w_+^q+\a g & \text{ in }\Omega, \\
u &=& 0  & \text{ on }\partial \Omega.\\
\end{array}
\right.
\end{equation}
It is clear that if $v$ is a fixed point of $T$, then $(u,v)$ solves the system \eqref{S}. Thus we just have to show that $T$ has a fixed point in $E$, this will be achieved in several steps.\vskip2mm
\noindent In what follows we denote by $C_1,C_2,...,$ any positive constants that depend only on the datum of the problem and that can be change from one line to next one.

\noindent {\bf Step I : $T$ is well defined.}
Let $w\in E$, by using Sobolev's inequality we conclude that $w\in L^{r^*}(\O)$. Since $q\s<r^*$, then $w_+^q\in L^\s(\O) \subset L^1(\O)$. Hence, by Theorem \ref{th0},   $u$ is well defined as the unique weak solution to the problem \eqref{P110} and $u\in W^{1, \theta}_0(\O)$ for all $\theta<\frac{N}{N-1}$. Taking into consideration the hypothesis on $g$, we reach that $(w^q_++\a g)\in L^{\s}(\O)$. Therefore, by Theorem \ref{key-2}, there exists $C_1$ such that
\begin{equation*}
 ||\n u||_{L^\b(\O)}\le C_1||w_+^q+\a g ||_{L^{\s}(\O)},\,\mbox{  for all } \b\leq \s^*=\frac{\s N}{N-\s}.
\end{equation*}
Using H\"older's inequality  and then by Sobolev's inequality, we get
\begin{equation}\label{eq1}
 ||\n u||_{L^\b(\O)}\le C_2\bigg(||\n w||^q_{L^{r}(\O)}+\a ||g||_{L^{\s}(\O)}\bigg),\,\mbox{  for all }  \b\leq \s^*.
\end{equation}
Since $\s^*>p$, then \eqref{eq1} holds with $\beta=p$, thus
\begin{equation}\label{eq9}
 ||\n u||_{L^p(\O)}\le C_2||\n w||^q_{L^{r}(\O)}+C_3||g||_{L^{\s}(\O)}.
\end{equation}
Since $f\in L^m (\O)$ and $|\n u|^p\in L^1(\O)$ then $|\n u|^p+\l f\in L^1(\O)$ and $v$ is the unique weak solution to the  problem \eqref{P100}. By applying again Theorem \ref{th0},   $v\in W^{1, \theta}_0(\O)$ for all $\theta<\frac{N}{N-1}$ and then $T$ is well defined.\\
\noindent{\textbf{Step II : $T(E)\subset E$}}.
First,   by  hypothesis $pm<\s^*$, we have $|\n u|^p\in L^{m}(\O)$. Second, by using the fact that $f\in L^{m}(\O)$ and by Theorem  \ref{key-2}, we obtain that for all $\theta\leq m^*=\frac{mN}{N-m}$,
\begin{equation}\label{ss}
\begin{array}{lll}
||\n v||_{L^{\theta}(\O)} & \le  & C_4\bigg\||\nabla u|^{p}+\l f\bigg\|_{L^{m}(\O)}\\
&\le & C_4\bigg(||\n u||^p_{L^{pm}(\O)}+\l ||f||_{L^{m}(\O)}\bigg).
\end{array}
\end{equation}
Recalling again that  $pm<\s^*$, then by using \eqref{eq1} with $\beta=pm$, we get
\begin{equation}\label{RR}
 ||\n u||^p_{L^{pm}(\O)}\le C_5\bigg(||\n w||^{pq}_{L^{r}(\O)}+\a^p ||g||^p_{L^{\s}(\O)}\bigg).
\end{equation}
Going back to \eqref{ss}, we conclude that
$$
||\n v||_{L^{\theta}(\O)}\le C_6\bigg(||\n w||^{pq}_{L^{r}(\O)}+\a^p ||g||^p_{L^{\s}(\O)}+\l ||f||_{L^{m}(\O)}\bigg).
$$
Since $r<\frac{mN}{N-m}$, by choosing $\theta=r$ in the previous inequality it holds that
$$
||\n v||_{L^{r}(\O)}\le C_7\bigg(||\n w||^{pq}_{L^{r}(\O)}+\a^p ||g||^p_{L^{\s}(\O)}+\l ||f||_{L^{m}(\O)}\bigg).
$$
Recall that $w\in E$, thus
$$
||\n v||_{L^{r}(\O)}\le C_8\bigg(\ell+\a^p ||g||^p_{L^{\s}(\O)}+\l ||f||_{L^{m}(\O)}\bigg).
$$
By choosing $\tilde{C}=C_8$ and taking into consideration the definition of $\ell$,
we conclude that\\
$||\n v||_{L^{r}(\O)}\le \ell^{\frac{1}{pq}}$. Thus $v\in E$ and then $T(E)\subset E$.\\

\noindent{\bf Step III : $T$ is a continuous and compact operator on $E$ endowed with the topology of $W^{1,1}_0(\O)$.}\\
%
$\bullet$ First, let us begin by proving the continuity of $T$. Assume that $\{w_n\}_n\subset E, w\in E$ are such that $w_n \to w \text{ in } W^{1,1}_0(\O)$ and define $v_n=T(w_n)$, $v=T(w)$. Then $(u_n,v_n)$  and $(u,v)$ satisfy

\begin{equation}\label{PP01}
\left\{
\begin{array}{rcll}
-\Delta u_n &= & w^q_{n+} +\a g & \text{ in }\Omega, \\
-\Delta u &= & w^q_++\a g & \text{ in }\Omega, \\
u_n&=& u= 0  & \text{ on }\partial \Omega,\\
\end{array}
\right.
\end{equation}
and
\begin{equation}\label{PP02}
\left\{
\begin{array}{rcll}
-\Delta v_n & = & |\nabla u_n|^{p}+\l f
& \text{ in }\Omega,\\
-\Delta v & = & |\nabla u|^{p}+\l f
& \text{ in }\Omega,\\
v_n &=& v=0  & \text{ on }\partial \Omega.
\end{array}
\right.
\end{equation}
Thanks to Sobolev's inequality in the space $W^{1,1}_0(\O)$, we have
$$
||w_n-w||_{L^{\frac{N}{N-1}}(\O)}\le S_1||w_n-w||_{W^{1,1}_0(\O)}\to 0\mbox{  as  }n\to \infty.
$$
Thus $w_n\to w$ strongly in $L^s(\O)$, for all $s\le \frac{N}{N-1}$. Since $\{w_n\}_n$ is bounded in the space $W^{1,r}_0(\O)$, then using Vitali's theorem it follows that $w_n\to w$ strongly in $L^a(\O)$ for all $a<r^*$. In particular $w^q_{n+}\to w^q_+$ strongly in $L^1(\O)$. Therefore thanks to Theorem \ref{th0} we conclude that $u_n\to u$ strongly in $W^{1,a}_0(\O)$ for all $a<\frac{N}{N-2}$. In particular we have
$$
||\n u_n-\n u||_{L^{1}(\O)}\to 0\mbox{  as  }n\to \infty.
$$
Now, by \eqref{RR} it follows that
\begin{equation}\label{RR10}
 ||\n u_n||^p_{L^{pm}(\O)}\le C\bigg(||\n w_n||^{pq}_{L^{r}(\O)}+\a^p ||g||^p_{L^{\s}(\O)}\bigg)\le C\bigg(\ell+\a^p ||g||^p_{L^{\s}(\O)}\bigg).
\end{equation}
Thus $\{u_n\}_n$ is bounded in $W^{1,pm}_0(\O)$. Using again Vitali's theorem we obtain that
$$
||\n u_n-\n u||_{L^{p}(\O)}\to 0\mbox{  as  }n\to \infty.
$$
Hence by Theorem \ref{th0} there results that $v_n\to v$ strongly in $W^{1,1}_0(\O)$ and then the continuity of $T$ follows.\\
$\bullet$ Second, we  show that $T$ is a compact operator.\\ Let $\{w_n\}_n\subset E$ be such that $||w_n||_{W^{1,1}_0(\O)}\le C$ and let $v_n=T(w_n)$. Since $\{w_n\}_n\subset E$,  $||\n w_n||_{L^{r}(\O)}\le C$ and then, up to a subsequence again denoted by $\{w_n\}_n$, we have
$$
w_n\rightharpoonup w \mbox{  weakly in } W^{1,r}_0(\O).
$$
By Rellich-Kondrachov's Theorem, it follows that $w_n\to w$ strongly in $L^a(\O)$ for all $a<r^*$. In particular, $w^q_{n+}\to w^q_+$ strongly in $L^1(\O)$.\\
Let $u$ to be the unique weak solution to the problem
\begin{equation}\label{CC01}
\left\{
\begin{array}{rcll}
-\Delta u &= & w^q_++\a g & \text{ in }\Omega, \\
u&=& 0  & \text{ on }\partial \Omega.
\end{array}
\right.
\end{equation}
By the result of Theorem \ref{th0}, we reach that $||\n u_n-\n u||_{L^{1}(\O)}\to 0$ as $n\to \infty.$ Now, setting $v=T(w)$ and following the same argument as in the proof of the continuity of $T$, we obtain that $v_n\to v$ strongly in $W^{1,1}_0(\O)$. Hence $T$ is compact. Therefore, by the Schauder Fixed Point Theorem we get the existence of $v\in E$ such that $T(v)=v$. It is clear that $v\in W^{1,r}_0(\O)$. Taking into consideration the hypothesis \eqref{mad}, it holds that $v\in W^{1,\rho}_0(\O)$ for all $\rho<\frac{mN}{N-m}$. Now going back to \eqref{P110}, we conclude that $u\in W^{1, \theta}_0(\O)$ for all $\theta<\frac{\s N}{N-\s}$. Thus the proof of Theorem \ref{syst1} follows. \cqd
\vskip3mm
\subsection{Proof of theorem \ref{smallness_Lambda-Alpha}}

Let $\alpha>0,\lambda >0$, $f\gneqq 0$,  $g\gneqq 0$. Assume that the system
\begin{equation}\label{S4}
\left\{
\begin{array}{rcll}
-\Delta u & = & v^q+\a g & \text{in }\Omega , \\
-\Delta v& = &|\nabla u|^{p}+\l f &\text{in }\Omega , \\
u=v&=& 0 & \text{on }\partial \Omega ,\\
u,v & \geq& 0 & \text{in }\Omega,
\end{array}%
\right.
\end{equation}
has a nonnegative solution $(u,v)$. Without loss of generality we can assume that $f,g\in L^\infty(\O)$.
\vskip2mm
\noindent In order to prove the existence of $\alpha^*$ and $\lambda^*$, we will distinguish two cases : $p,q >1$ and $p>1, q= 1$.
\vskip2mm
\noindent  In what follows we will use Young's inequality under this form
$$\forall (a,b)\in(0,+\infty)^2 \text{ and } \forall  s \in(1,+\infty), \quad
ab\leq a^s+C_sb^{s'} \text{ where } C_s= \frac{s-1}{s^{s'}} \text{ and } s'=\frac{s}{s-1}.$$
$\bullet$ \noindent \textbf{Case $p,q >1$.}\\
Let $\phi\in \mathcal{C}_0^{\infty}(\Omega)$ such $\phi \gneqq 0$.
Multiplying both equations of (\ref{S4}) by $\phi$ and integrating over $\Omega$, we obtain
\begin{equation}\label{pour-controler-alpha-g-1}
\int_\Omega \left(v^q +\alpha g\right)\phi= \int_\Omega \nabla u\nabla\phi
\end{equation}
and
\begin{equation}\label{pour-controler}
\int_\Omega \left(|\nabla u|^p +\lambda f\right)\phi= -\int_\Omega \phi\Delta v = \int_\Omega  v(-\Delta \phi).
\end{equation}
\noindent \textbf{ First we will determine } $\alpha ^*$.\\
Thanks to  Young's inequality,  (\ref{pour-controler-alpha-g-1}) yields
\begin{equation}\label{pour-controler-alpha-g-2}
\int_\Omega \left(v^q +\alpha g\right)\phi
\leq
\int_\Omega \phi |\nabla u|^p  + C_p\int_\Omega \phi^{1-p'} |\nabla \phi|^{p'}
\end{equation}
By  the nonnegativity of each term of (\ref{pour-controler})
we get
\begin{equation}\label{pour-controler-alpha-g-3}
\int_\Omega|\nabla u|^p\phi\leq \int_\Omega  v(-\Delta \phi).
\end{equation}
 Therefore Young's inequality yields
\begin{equation}\label{pour-controler-alpha-g-4}
\int_\Omega|\nabla u|^p\phi
\leq
\int_{\Omega}v^q\phi+C_q\int_{\Omega} \phi^{1-q'}|\Delta \phi|^{q^{\prime}}.
\end{equation}
Let us set
$$F(\phi)=C_p\int_\Omega \phi^{1-p'} |\nabla \phi|^{p'} + C_q\int_{\Omega} \phi^{1-q'}|\Delta \phi|^{q^{\prime}}.$$
From (\ref{pour-controler-alpha-g-2}) and (\ref{pour-controler-alpha-g-4}), we deduce
\begin{equation}\label{pour-controler-alpha-g-5}
\alpha\int_{\Omega}g \phi
\leq F(\phi)
\end{equation}
and then $\a\le \a^*$ with
\begin{equation}\label{alpha-etoile}
\alpha ^*=\inf\left\{F(\phi)\;;\; 0\leq\phi\in \mathcal{C}^{\infty}_0(\Omega) \text{ and } \displaystyle\int_{\Omega}g\phi=1 \right\}.
\end{equation}
\vskip2mm
\noindent It is clear that  if $\alpha > \alpha ^*$, the system (\ref{S}) has no nonnegative solution.
\vskip2mm
\noindent \textbf{ Now we will determine $\lambda ^*$.}\\ For this, we will proceed in the same way as before.\\
By applying Young's inequality, (\ref{pour-controler}) implies
\begin{equation}\label{pour-controler-lambda-f-1}
\int_\Omega \left(|\nabla u|^p +\lambda f\right)\phi
\leq
 \int_\Omega \phi v^q  + C_q\int_\Omega \phi^{1-q'} |\Delta \phi|^{q'}
\end{equation}
Thanks to the nonnegativity of each term, (\ref{pour-controler-alpha-g-1}) yields
\begin{equation}\label{pour-controler-lambda-f-3}
\int_\Omega v^q\phi\leq\int_\Omega \nabla u \nabla \phi
\end{equation}
and then
\begin{equation}\label{pour-controler-lambda-f-4}
\int_\Omega v^q\phi\leq
\int_\Omega \phi |\nabla u|^p  + C_p\int_\Omega \phi^{1-p'} |\nabla \phi|^{p'}
\end{equation}
 Therefore (\ref{pour-controler-lambda-f-1}) and (\ref{pour-controler-lambda-f-4}) imply
 \begin{equation}\label{pour-controler-lambda-f-5}
\lambda \int_\Omega  f\phi
\leq
F(\phi).
\end{equation}
and then
\begin{equation}\label{lambda-etoile}
\lambda ^*=\inf\left\{F(\phi)\;;\; 0\leq\phi\in \mathcal{C}^{\infty}_0(\Omega) \text{ and } \displaystyle\int_{\Omega}f\phi=1 \right\}.
\end{equation}
%
%
It is clear that  if $\lambda > \lambda ^*$, the system (\ref{S}) has no nonnegative solution.
%
%
%
\vskip4mm
To complete the proof in this part we have to show that $\a^*, \l^*>0$.
\vskip2mm
\noindent Notice that for $\phi$ as above, we have
$$
F(\phi)\ge C_p\int_\Omega \phi^{1-p'} |\nabla \phi|^{p'}=C_p(p')^{p'}\int_\Omega|\nabla \phi^{\frac{1}{p'}}|^{p'}.
$$
Without loss of generality we can assume that $p'<N$, thus, using Sobolev's inequality we reach that
$$
F(\phi)\ge C_p\bigg(\int_\Omega |\phi|^{\frac{N}{N-p'}}\bigg)^{\frac{N-p'}{N}}.
$$
Hence, using H\"older's inequality, we conclude that
$$
\dfrac{F(\phi)}{\dyle\io g\phi}\ge \dfrac{C_p}{||g||_{L^{\frac{N}{p'}}(\O)}}.
$$
Therefore we conclude that $\a^*\ge \dfrac{C_p}{||g||_{L^{\frac{N}{p'}}(\O)}}.$ In the same way we reach that $\l^*>0$.
%
%
\vskip3mm
\noindent $\bullet$ \noindent\textbf{ Case $q=1$ and $p>1$.}\\
%
%
Let $\varphi\in \mathcal{C}^\infty_0(\O)$ be such that $\varphi\gneqq 0$ and define $\phi$ to be the unique solution to the following problem
\begin{equation}\label{inter}
\left\{
\begin{array}{rcll}
-\Delta \phi & = & \varphi & \text{in }\Omega , \\
\phi&=& 0 & \text{on }\partial \Omega.
\end{array}%
\right.
\end{equation}
It is clear that $\phi\in \mathcal{C}^2(\bar{\O})$ and $-\Delta \phi \gneqq 0$. Then multiplying the first equation in \eqref{S} by $-\Delta \phi$, the second equation  by $\phi$  and integrating to obtain
\begin{equation}\label{pour-controler-alpha-g-1-1}
\int_{\Omega}v(-\Delta \phi)+\alpha\int_{\Omega}g(-\Delta\phi)=\int_{\Omega}(-\Delta u)(-\Delta \phi)=\int_{\Omega}\nabla u.\nabla (-\Delta \phi)\\
\end{equation}
\textbf{  Determination of $\alpha^*$. } \\
Thanks to Young's inequality, we have
\begin{equation}\label{pour-controler-alpha-g-2-1}
\int_{\Omega}v(-\Delta \phi)+\alpha\int_{\Omega}g(-\Delta\phi)
\leq
\int_{\Omega}|\nabla u|^p\phi+C_p\int_{\Omega}\phi^{1-p^{\prime}}|\nabla (-\Delta \phi)|^{p^{\prime}}.
\end{equation}
But
$$\int_{\Omega}|\nabla u|^p\phi\leq \int_{\Omega}|\nabla u|^p\phi +\lambda\int_{\Omega}f\phi=\int_{\Omega}(-\Delta v)\phi=\int_{\Omega}v(-\Delta\phi)$$
Therefore
\begin{equation}\label{pour-controler-alpha-g-3-1}
\int_{\Omega}v(-\Delta \phi)+\alpha\int_{\Omega}g(-\Delta\phi)
\leq
\int_{\Omega}(-\Delta v)\phi + C_p\int_{\Omega}\phi^{1-p^{\prime}}|\nabla (-\Delta \phi)|^{p^{\prime}}.
\end{equation}
Thus
\begin{equation}\label{pour-controler-alpha-g-4-1}
\alpha\int_{\Omega}g(-\Delta\phi)
\leq C_p\int_{\Omega}\phi^{1-p^{\prime}}|\nabla (-\Delta \phi)|^{p^{\prime}}.
 \end{equation}
Going back to the definition of $\phi$ and  setting $G(\varphi)= C_p\displaystyle\int_{\Omega}\phi^{1-p^{\prime}}|\nabla \varphi|^{p^{\prime}}$, it holds that $\a\le \a^*$
 with
\begin{equation}\label{pour-controler-alpha-g-5-1}
\alpha^*=\inf\left\{G(\varphi)\; | \; \varphi\in \mathcal{C}^\infty_0(\O)\;,\; \varphi\gneqq 0 \text{ and }\displaystyle \int_\Omega g\varphi=1\right\}.
\end{equation}
\noindent\textbf{  Determination of $\lambda^*$. }\\
In the same way, we obtain

\begin{equation}\label{pour-controler-lambda-f-1-1}
\int_{\Omega}|\nabla u|^p\phi+\lambda\int_{\Omega}f\phi=\int_{\Omega}(-\Delta v) \phi=\int_{\Omega}v(-\Delta \phi)
\end{equation}
and
\begin{equation}\label{pour-controler-lambda-f-2-1}
\int_{\Omega}v(-\Delta \phi)\leq \int_{\Omega}|\nabla u|^p\phi+ C_p\int_{\Omega}\phi^{1-p^{\prime}}|\nabla (-\Delta\phi)|^{p^{\prime}}.
\end{equation}
%
Thus
\begin{equation*}\label{pour-controler-lambda-f-2-111}
\lambda\int_{\Omega}f\phi \leq  G(\varphi),
\end{equation*}
and then $\l\le \l^*$ with
\begin{equation}\label{pour-controler-lambda-f-3-1}
\lambda^*=\inf\left\{G(\varphi)\;|\; \int_\Omega f \phi=1\right\}.
\end{equation}
\vskip3mm

As in the first case we will show that $\a^*,\l^*>0$. But the situation here is more complicated and needs some fine estimates.\\
Recall that $\phi$ is the unique solution of problem \eqref{inter}. Fix $1<\theta<\min\{\frac{N}{(N-3(p'-1))_+}, p'\}$, then $\frac{(p'-1)\theta}{p'-\theta}\le \frac{\theta N}{N-3\theta}$.
But $\varphi\in \mathcal{C}^{\infty}_0(\O)$ implies $\varphi\in W^{1,p'}_0(\O)$. Thus, taking into consideration the definition of $\phi$ and $\theta$, and by Theorem \ref{key-2}, it follows that
$$
||\phi||_{L^{\frac{(p'-1)\theta}{p'-1}}(\O)}\le C||\varphi||_{L^{\theta^*}(\O)}\le C||\n \varphi||_{L^{\theta}(\O)}.
$$
Now, using H\"older's inequality, it follows that
\begin{eqnarray*}
\int_{\Omega}|\nabla \varphi|^{\theta} &= & \dyle\int_{\Omega}|\nabla \varphi|^{\theta}\phi^{\frac{(1-p')\theta}{p'}}\; \phi^{\frac{(p'-1)\theta}{p'}}\le \dyle \bigg(\int_{\Omega}\phi^{1-p^{\prime}}|\nabla \varphi|^{p^{\prime}}\bigg)^{\frac{\theta}{p'}}
\bigg(\int_{\Omega}|\varphi|^{\frac{(p'-1)\theta}{p'-\theta}}\bigg)^{\frac{p'-\theta}{p'}}\\
&\le & \dyle C\bigg(\int_{\Omega}\phi^{1-p^{\prime}}|\nabla \varphi|^{p^{\prime}}\bigg)^{\frac{\theta}{p'}}
\bigg(\int_{\Omega}|\varphi|^{\theta^*}\bigg)^{\frac{(p'-1)\theta}{p'\theta^*}}
\end{eqnarray*}
Thus, using Sobolev's inequality it holds that
$$
\bigg(\int_{\Omega}|\nabla \varphi|^{\theta} \bigg)^{\frac{1}{\theta}}\le C\int_{\Omega}\phi^{1-p^{\prime}}|\nabla \varphi|^{p^{\prime}}=CG(\varphi).
$$
Hence
$$
\dyle \dfrac{G(\varphi)}{\dyle\io g\varphi}\ge \dyle \dfrac{C(p,\O)\bigg(\dyle\int_{\Omega}|\nabla \varphi|^{\theta} \bigg)^{\frac{1}{\theta}}}{\dyle\io g\varphi}\ge \dfrac{C(p,\O)}{||g||_{L^{\frac{\theta N}{(\theta-1)N+\theta}}(\O)}}.
$$
Thus $\a^*>0$.

\

Now, to deal with $\l^*$, we again apply H\"older's inequality to obtain
$$
\io f\phi\le \bigg(\io \phi^{\frac{\theta N}{N-3\theta}}\bigg)^{\frac{N-3\theta}{\theta N}} \bigg(\io f^{\frac{\theta N}{(\theta-1)N+3\theta}}\bigg)^{\frac{(\theta-1)N+3\theta}{\theta N}}.
$$
Taking into consideration that $\phi$ solves \eqref{inter}, we conclude that
$$
\bigg(\io \phi^{\frac{\theta N}{N-3\theta}}\bigg)^{\frac{N-3\theta}{\theta N}}\le C\bigg(\io \varphi^{\theta^*}\bigg)^{\frac{1}{\theta^*}}.
$$
Thus
$$
\dyle \dfrac{G(\varphi)}{\dyle\io f\phi}\ge \dyle \dfrac{C\bigg(\dyle\int_{\Omega}|\nabla \varphi|^{\theta} \bigg)^{\frac{1}{\theta}}}{\bigg(\dyle\io \varphi^{\theta^*}\bigg)^{\frac{1}{\theta^*}}||f||_{L^{\frac{\theta N}{N-3\theta}}(\O)}}\ge \dfrac{C(N,\O)}{||f||_{L^{\frac{\theta N}{N-3\theta}}(\O)}},
$$
then $\l^*>0$. \cqd

\begin{remark}
The existence result in Theorem \ref{syst1} is obtained under the hypothesis that $pq>1$ and smallness condition on $(\a,\l)$. It will be interesting to show that nonexistence result, for large values of $(\a,\l)$, obtained in Theorem \ref{smallness_Lambda-Alpha} holds under the condition $pq>1$.
\end{remark}
\vskip3mm
%
%
\subsection{Proof of Theorem \ref{bih}}
%

In this case we take $q=1$ and $g=0$ in Theorem \ref{syst1}. Since $f\in L^s(\O)$ with $m>\max\{1,\frac{N}{3p'}\}$, we can fix $m_0$ such that $m_0\le m$ and $\frac{N}{3p'}<m_0<N$.

If $m_0\ge \frac{N}{2}$, we take $\s_0=N-\e$ with $\e$ small enough, however, if $\frac{N}{3p'}<m_0<\frac{N}{2}$, we choose $\s_0=\frac{3pm_0}{p+2}$. Hence in all cases, the condition \eqref{condi0} holds and then by Theorem \ref{syst1}, we get the existence of $(u,v)\in W^{1,\rho}_0(\O)\times W^{1, \theta}_0(\O)$, for all $\rho<\frac{m_0N}{N-m_0}, \theta<\frac{\s_0 N}{N-\s_0}$, such that $(u,v)$ solves the system
\begin{equation*}
\left\{
\begin{array}{rcll}
-\Delta u & = & v & \text{in }\Omega , \\
-\Delta v& = &|\nabla u|^{p}+\l f &\text{in }\Omega , \\
u=v&=& 0 & \text{on }\partial \Omega ,\\
u,v& \geq & 0 & \text{in }\Omega.
\end{array}%
\right.  \label{P0}
\end{equation*}
It is clear that:
$$
u=\D u=0 \quad \hbox{ on } \p\O \quad \hbox{ and }\quad
\D^2 u=-\D v=|\n u|^p+\l f\quad \hbox{ in } \O.
$$
Hence $u$ solves \eqref{PBi}. Now using the regularity results in \cite{GG}, we reach that $u\in \mathcal{C}^s(\Omega)$ for all $s<4-\frac{N}{m}$, in particular $u\in \mathcal{C}^s(\Omega)$ for all $s<2$
if $m\ge \frac{N}{2}$ and $u\in W^{4,m}(\O)$ if $m<\frac{N}{2}$. Hence we conclude. \cqd

\begin{remark}
Taking into consideration the nonlinear nature of the above arguments used in the proof of Theorem \ref{syst1}, then we can extend our existence result for more general operator like the $m$-laplacian operator where $\Delta_m(u): =\text{div}(|\nabla u|^{m-2}\nabla u)$. More precisely if we consider the system
\begin{equation}\label{Plapl}
\left\{
\begin{array}{rcll}
-\Delta_{m_1} u & = & v^q+\a g & \text{in }\Omega , \\
-\Delta_{m_2} v& = &|\nabla u|^{p}+\l f &\text{in }\Omega , \\
u=v&=& 0 & \text{on }\partial \Omega ,\\
u,v & \geq & 0 & \text{in }\Omega,
\end{array}%
\right.
\end{equation}
with $m_1, m_2>2-\frac{1}{N}$, then if $f\in L^\theta(\O)$ and $g\in L^\s(\Omega)$, where  $(\theta,\s)$ satisfies
\begin{equation}\label{condi011}
\left\{
\begin{array}{rcll}
\theta,\s &\in & (1, N),\\ \\
p\theta & < & \dfrac{\s N(m_1-1)}{N-\s},\\ \\
\frac{q\s N}{N+q\s} &<& \dfrac{\theta N(m_2-1)}{N-\theta},
\end{array}%
\right.
\end{equation}
we can show the existence of a bounded set $\overline{\Pi}\subset \re^+\times \re^+$ defined by
$$
\overline{\Pi}\equiv \bigg\{(\l, \a)\in (\re^+)^2: \l||f||_{L^{\theta}(\O)}+\a^p||g||_{L^{\s}(\O)}^p\le \overline{\L}\bigg\},
$$
such that for all $(\l, \a)\in \overline{\Pi}$, the system \eqref{Plapl} has a positive solution $(u,v)$ with $(u,v)\in W^{1,\rho}_0(\O)\times W^{1, r}_0(\O)$ for all $\rho<\frac{\s N(m_1-1)}{N-\s}$ and $r<\frac{\theta N(m_2-1)}{N-\theta}$.
\end{remark}
\vskip3mm
\subsection{Some nonexistence results}
%
%
%

In this subsection we will show that the regularity conditions \eqref{condi0} are optimal for existence at least in some explicit cases. More precisely we will analyze the problem \eqref{PBi}, that is the case where $q=1$ and $g\equiv 0$.
\vskip3mm
Let us begin by proving the following result.
\begin{Theorem}\label{Non1}
Assume that $0\le f\in L^m(\O)$ with $m\ge 1$ and $p>1$. Define
\begin{equation}\label{RRR}
\L(f)\equiv \inf_{\phi\in \mathcal{C}^\infty_0(\O)}\dfrac{\dyle\io \dfrac{|\n \big(|\phi|^{p'-2}\phi (-\D \phi)\big)|^{p'}}{|\phi|^{\frac{p}{(p-1)^2}} }dx}{\dyle\io f|\phi|^{p'}dx}.
\end{equation}
Assume that problem \eqref{PBi} has positive solution for some $\l>0$, then $\L(f)>0$.
\end{Theorem}
%
\textbf{\em Proof :} We follow closely the arguments used in \cite{BPihp}.
Suppose that problem \eqref{PBi} has positive solution $u$. Setting $v=-\D u$, then $v$ solves
$$
\left\{
\begin{array}{rcll}
-\Delta v &=& |\n u|^p+ \l f & \mbox{ in }\Omega, \\ v &=& 0 & \mbox{  on } \p\Omega
\end{array}%
\right.
$$
Thus $v>0$, hence $-\D u>0$ in $\O$. Let $\phi\in \mathcal{C}^\infty_0(\O)$, using $|\phi|^{p'}$ as a test function in \eqref{PBi}, it follows that
\begin{equation}\label{test1}
\io (-\D u)(-\D |\phi|^{p'})dx =\io |\n u|^p|\phi|^{p'} dx +\l \io f |\phi|^{p'}dx.
\end{equation}
Since $p'>1$, then using Kato's inequality we reach that
$$
-\D |\phi|^{p'}\le p'|\phi|^{p'-2}\phi (-\D \phi).
$$
Since $-\D u\ge 0$, using Young's inequality, it follows that
\begin{eqnarray*}
\io (-\D u)(-\D |\phi|^{p'})dx &\le & p'\io |\phi|^{p'-2}\phi (-\D \phi)(-\D u)dx\\
&\le & p'\io \n u\n \big(|\phi|^{p'-2}\phi (-\D \phi)\big)dx\\
&\le & \e \io |\n u|^p|\phi|^{p'} dx +C(\e)\io \dfrac{|\n \big(|\phi|^{p'-2}\phi (-\D \phi)\big)|^{p'}}{|\phi|^{\frac{p}{(p-1)^2}}}dx.
\end{eqnarray*}
Choosing $p'\e\le 1$ and going back to \eqref{test1}, we conclude that
$$
\l \io f |\phi|^{p'}dx\le C(\e,p)\io \dfrac{|\n \big(|\phi|^{p'-2}\phi (-\D \phi)\big)|^{p'}}{|\phi|^{\frac{p}{(p-1)^2}}}dx.
$$
Thus
\begin{equation}\label{yyy}
\l\le C(\e,p) \dfrac{\dyle\io \dfrac{|\n \big(|\phi|^{p'-2}\phi (-\D \phi)\big)|^{p'}}{|\phi|^{\frac{p}{(p-1)^2}}}dx}{\dyle\io f|\phi|^{p'}dx}
\end{equation}
and the result follows. \cqd

\begin{remark}
Let $\phi\in \mathcal{C}^\infty_0(\O)$ and define
$$
Q(\phi)=\dfrac{\dyle\io \dfrac{|\n \big(|\phi|^{p'-2}\phi (-\D \phi)\big)|^{p'}}{|\phi|^{\frac{p}{(p-1)^2}}}dx}{\dyle\io f|\phi|^{p'}dx},
$$
then $Q(\a \phi)=Q(\phi)$ for all $\a\in \re^*$.
\end{remark}

As a direct consequence of Theorem \ref{Non1} we can show that problem \eqref{PBi} has non solution if $\l$ is large. This follows directly by estimate \eqref{yyy}.
\vskip3mm

We are now able to prove the next nonexistence result.
\begin{Theorem}\label{Non1001}
Assume that $1\le m < \max\{1,\frac{N}{3p'}\}$ where $p>1$, then there exists $f\in L^m(\O)$ with $f\gneq 0$ such that problem \eqref{PBi} has non positive solution for all $\l>0$.
\end{Theorem}
%
\textbf{\em Proof :} Without loss of generality we can assume that $N>p'$. Assume that $\O=B_1(0)\subset \ren$ and fix $m<\frac{N}{3p'}$. Define the function $f(x)=\dfrac{1}{|x|^{3p'+\e}}$, then we can choose $\e$ small enough such that $f\in L^m(B_1(0))$. To show the non existence result we will prove that $\L(f)=0$ and then by the previous Theorem we will conclude.

Define the function $\phi$ be
$$
\phi(x)=
\left\{
\begin{array}{lll}
\dfrac{1}{|x|^\theta}& \mbox{  if  }& |x|\le \frac 14\\
(1-|x|)^\g & \mbox{  if  }& \frac 12\le |x|\le 1,
\end{array}
\right.
$$
where $\theta=\frac{N-(3p'+\e)}{p'}$, $\g>\frac{3p'-1}{p'}$, and $\phi\in \mathcal{C}^2(B_1(0)\setminus\{0\})$ with $\phi>0$ in $B_1(0)$.

Now by a direct computation we obtain that
$$
\int_{B_1(0)} \dfrac{|\n \big(|\phi|^{p'-2}\phi (-\D \phi)\big)|^{p'}}{|\phi|^{\frac{p}{(p-1)^2}}}dx=\int_{B_{\frac 14}(0)}+\int_{\frac 14<|x|<\frac 12}
+\int_{\frac 12<|x|<1}.
$$
Since $\phi>0$ for $\frac 14<|x|<\frac 12$, then $\int_{\frac 14<|x|<\frac 12}=C<\infty$.

By a direct computations we reach that
$$
\int_{B_{\frac 14}(0)}\dfrac{|\n \big(|\phi|^{p'-2}\phi (-\D \phi)\big)|^{p'}}{|\phi|^{\frac{p}{(p-1)^2}}}dx\le C\int_{B_{\frac 14}(0)}\dfrac{1}{|x|^{p'(\theta+3)}}dx.
$$
Since $\theta=\frac{N-(3p'+\e)}{p'}$, then $p'(\theta+3)=N-\e<N$. Thus
$$
\int_{B_{\frac 14}(0)}\dfrac{|\n \big(|\phi|^{p'-2}\phi (-\D \phi)\big)|^{p'}}{|\phi|^{\frac{p}{(p-1)^2}}}dx\le C.
$$
We deal now with the integral $\dyle\int_{\frac 12<|x|<1}$. In this case we have
$$
\int_{\frac 12<|x|<1}\dfrac{|\n \big(|\phi|^{p'-2}\phi (-\D \phi)\big)|^{p'}}{|\phi|^{\frac{p}{(p-1)^2}}}dx\le C\int_{\frac 12<|x|<1}(1-|x|)^{p'(\g-3)}dx.
$$
Since $\g>\frac{3p'-1}{p'}$, then $p'(\g-3)>-1$. Thus
$$
\int_{\frac 12<|x|<1}\dfrac{|\n \big(|\phi|^{p'-2}\phi (-\D \phi)\big)|^{p'}}{|\phi|^{\frac{p}{(p-1)^2}}}dx\le C.
$$
As a conclusion we have proved that
$$
\int_{B_1(0)} \dfrac{|\n \big(|\phi|^{p'-2}\phi (-\D \phi)\big)|^{p'}}{|\phi|^{\frac{p}{(p-1)^2}}}dx<\infty.
$$
On the other hand we have $\dyle\io f|\phi|^{p'}dx=+\infty$. Hence $Q(\a \phi)=0$. Now using an approximation argument we obtain that $\L(f)=0$ and therefore the result follows.\cqd

\begin{remark}
Combining the existence result in Theorem \eqref{bih} and the result of Theorem \ref{Non1} we obtain that if $0\lneqq f\in L^m(\O)$ where  $m >\max\{1,\frac{N}{3p'}\}$, then $\L(f)>0$.

It will be nice if we can prove directly that if $f\in L^m(\O)$ with $m >\max\{1,\frac{N}{3p'}\}$, then $\L(f)>0$.
\end{remark}

\noindent {\bf Acknowledgments :}\\
$\bullet$ The authors would like to thank Prof. Michel  PIERRE for his helpful suggestions and fruitful discussions during the preparation of this work.
\vskip3mm
\noindent $\bullet$ Part of this work was realized while the first author was visiting the Institute Elie Cartan, Universit\'{e} Lorraine. He would like to thanks the institute for the warm hospitality.


\begin{thebibliography}{99}

\bibitem{ADP} B. Abdellaoui, A. Dall'Aglio,  and I. Peral,  {\it Some
	remarks on elliptic problems with critical growth in the
	gradient}. J. Differential Equations {\bf 222} (2006), 21-62.
	
\bibitem{AdamsPierre}  D. R. Adams, M. Pierre, {\em Capacitary strong type estimates in semilinear problems}. Ann. Inst. Fourier (Grenoble) {\bf 41} (1991), no. 1, 117--135.

\bibitem{Alp} N.E. Alaa, M. Pierre, {\it Weak solutions of some quasilinear elliptic equations with data measures}.
SIAM J. Math. Anal., {\bf 24},  (1993), 23-35.

\bibitem{AtB} A. Attar, R. Bentifour, {\em Nonlinear elliptic system involving gradient term and reaction potential }. Electronic Journal of Differential Equations, {\bf 2017} (2017), no. 113, 1-10.

\bibitem{BPihp} P. Baras, M. Pierre, {\it Crit{è}re d'existence des solutions positives pour des équations semi-linéaires
non monotones}, Ann. I.H.P. {\bf 2}, (3) (1985), 185-212.

\bibitem{BPaif} P. Baras, M. Pierre, {\it Singularités éliminables pour des équations semi-linéaires,}
Ann. Inst. Fourier {\bf 34}, (1984), no. 1, 185-206.

\bibitem{BarlPorre}
G. Barles, A. Porretta, {\it Uniqueness for unbounded solutions to stationary viscous Hamilton-Jacobi equations.}
 Ann. Sc. Norm. Super. Pisa Cl. Sci. {\bf 5} (2006), no. 1, 107-136.



\bibitem{BocCroce}L. Boccardo, G. Croce, {\it Elliptic Partial Differential Equations : existence and regularity of distributional solutions}. Studies in Mathematics \textbf{55} (2014), De Gruyter.

\bibitem{BGM} L. Boccardo, T. Gallou\"{e}t, F. Murat, {\it A unified
representation of two existence results for problems with natural
growth}, Research Notes in Mathematics {\bf 296} (1993),
127-137.

\bibitem{BMPu}  L. Boccardo, F. Murat, J.-P. Puel, {\it Existence des
solutions non born\'ees pour certains \'equations
quasi-lin\'eaires},  Portugal Math. {\bf  41} (1982), 507-534.


\bibitem{BOPorr} L. Boccardo, L. Orsina, A. Porretta, {\it  Existence of finite energy solutions for elliptic systems with $L^1-$ value nonlinearities} Mathematical Models and Methods in Applied Sciences {\bf 18}, (2008), no.5, 669-687.


\bibitem{BOPu} L. Boccardo, L. Orsina, J.-P. Puel, {\it A quasilinear elliptic system with natural growth terms}.  Annali di Matematica.  {\bf 194}, (2015), no.3, 1733-1750.

\bibitem{Bre} H. Br\'ezis, {\it  Functional analysis, Sobolev spaces and partial differential equations.}  Universitext. Springer, New York, 2011.


\bibitem{ChoCho} K. Cho, H.J. Choe, {\it Nonlinear degenerate elliptic partial differential equations with critical growth conditions on the
gradient}, Proc. A.M.S. {\bf  123}, (1995), no 12, 3789-3796.

\bibitem{C} S. Clain, J. Rappaz, M. Swierkosz, R. Touzani, {\it Numerical modeling of induction heating for two dimentional geometrie,}
 Math. Models Methods Appl. Sci. {\bf 3}, (1993), no. 6, 805-822.



\bibitem{Diez} J. I. Díaz, M. Lazzo, P.G. Schmidt, {\it Large solutions
for a system of elliptic equation arising from fluid dynamics}. Siam Journal on Mathematical Analysis, {\bf 37} (2005), 490-513.

\bibitem{EP} C. Escudero, I. Peral, {\it Some fourth order nonlinear elliptic problems related to epitaxial growth}, J.
Differential Equations, {\bf 254} (2013), 2515-2531.

\bibitem{FM1} V. Ferone, F. Murat, {\it Quasilinear problems having
quadratic growth in the gradient : an existence result when the
source term is small}, Equations aux d\'eriv\'ees partielles et
applications,  497-515, Gauthier-Villars, Ed. Sci. M\'ed.
Elsevier, Paris, 1998.




\bibitem{GG} F. Gazzola, H.C. Grunau, G. Sweers, {\it Polyharmonic boundary value problems. Positivity preserving and nonlinear higher order elliptic equations in bounded domains.} Lecture Notes in Mathematics, 1991. Springer-Verlag, Berlin, 2010.


\bibitem{GMP} N. Grenon, F. Murat, A. Porretta, {\it Existence and a priori
estimate for elliptic problems with subquadratic gradient
dependent terms}, C. R. Acad. Sci. Paris, Ser. I {\bf 342}
(2006), 23-28.

\bibitem{GreTrom} N. Grenon, C. Trombetti, {\it Existence results for a class of nonlinear elliptic problems with $p$-growth in the gradient},
Nonlinear Anal. {\bf  52} (2003), no.3, 931-942.

\bibitem{HMV} K Hansson, V.G. Maz'ya, I.E. Verbitsky, {\it Criteria of solvability for multidimensional
Riccati equations}, Ark. Mat., {\bf 37}, (1999), 87-120.

\bibitem{LL} J. Leray and J.-L. Lions, {\it Quelques r\'esultats
de Vi\v{s}ik sur les probl\`{e}mes elliptiques non lin\'eaires par
les m\'ethodes de Minty-Browder}, Bull. Soc. Math. France {\bf  93},
(1965), 97-107.

\bibitem{PLL}  P.-L. Lions,{\it Generalized solutions of Hamilton-Jacobi
Equations}, Pitman Res. Notes Math. {\bf 62} (1982).

\bibitem{Por} A. Porretta, {\it Nonlinear equations with natural
growth terms and measure data}, 2002-Fez Conference on Partial
Differential Equations, Electron. J. Diff. Eq. Conf. {\bf 09} (2002),
183-202.


\bibitem{Phuc} N.C. Phuc, {\it  Morrey global bounds and qusilinear Riccarti type equation bellow the natural exponent}, J. math. Pures Appl. {\bf 102}, (2014), 99-123.


\bibitem{St} G. Stampacchia, {\it Le probl\`eme de Dirichlet pour
les \'equations elliptiques du second ordre \`a    coefficients discontinus},  Ann. Inst. Fourier (Grenoble), {\bf 15}(1965), 189-258.


\end{thebibliography}
\end{document}